\documentclass{article}

\usepackage[scaled]{helvet}
 
\usepackage[T1]{fontenc}

\usepackage{graphicx,pst-text}
\usepackage{pstricks,pst-node,amssymb,relsize,amsmath,wrapfig,amsmath,amsthm,url}
\newcommand{\jump}[1]{\left[\kern-3pt\left[{#1}\right] \kern-3pt\right]}
\newcommand{\jumpa}[1]{\left[\kern-2pt\left[{#1}\right] \kern-2pt\right]}
\DeclareMathAlphabet{\sans}{OT1}{cmss}{bx}{n}

\def\L{\mathcal{L}}

\def\SW{SW equations}
\def\SG{SG equations}
\newcommand{\R}{\mathbb{R}}

\usepackage{chngcntr}
\counterwithout{figure}{section}

\newtheoremstyle{example}{\topsep}{\topsep}%
     {}
     {}
     {\bfseries}
     {}
     {\newline}
     {\thmname{#1}\thmnumber{ #2}\thmnote{ #3}}

   \theoremstyle{example}
   \newtheorem{example}{Example}[section]
\newtheorem{remark}{Remark}[section]
\newtheorem{lemma}{Lemma}[section]

\def\ds{\displaystyle}

\begin{document}

\title{ On discrete analogues of potential vorticity for\\ variational shallow water systems }

\author{ Elizabeth L.\ Mansfield\\  School of Mathematics, Statistics and Actuarial Science,\\ University of Kent, Canterbury, CT2 7FS, UK}
\date{March 7, 2016}
\maketitle

\begin{abstract}
We outline how discrete analogues of the conservation of potential vorticity may be achieved in Finite Element numerical schemes for 
a variational system which has the particle relabelling symmetry, typically shallow water equations.

We show that the discrete analogue of the conservation law for potential vorticity converges to the smooth law for
potential vorticity, and moreover, for a strong solution, is the weak version of the potential vorticity law. This result rests
on recent results by the author with T.\ Pryer concerning discrete analogues of conservation laws in Finite Element variational problems, together with 
an observation by P.\ Hydon concerning how the conservation of potential vorticity in smooth systems arises as a consequence of the linear momenta.


The purpose of this paper is to provide all the necessary information for the implementation of the schemes and the necessary numerical tests. A brief tutorial on Noether's theorem is included
to demonstrate the origin of the laws and to demonstrate that the numerical method follows the same basic principle, which is that the law follows directly from the Lie group invariance of the Lagrangian. 
\end{abstract}

\section{Introduction}
The UK Met office famously failed to predict a hurricane in 1987. The author was told that\footnote{The author is grateful to Ian Roulstone, formerly of the Met Office, for relating the story to her.}, in retrospect, this was felt to be the failure of the numerical code to preserve  the physically important
quantity, {potential vorticity} (as well as incomplete data).

More recently, in 2012, a post-tropical cyclone hit the United States, near Brigantine, New Jersey. Several early forecasts of this cyclone stated that the weather system would disappear over the Atlantic Ocean. However, researchers at the European Center for Medium-Range Weather Forecasts (ECMWF) had a different opinion: they possessed a model which predicted the landfall of the storm around New Jersey. Their successful forecast derived from the fact that their model depicted
accurately the evolution of potential vorticity  \cite{RN13}. Thus, potential vorticity plays an important role in weather forecasting, and for this reason, it is interesting to investigate its inclusion in a numerical model.

Potential vorticity arises as a Noether conservation law under the ``particle relabelling" symmetry, that is, invariance of the Lagrangian under area\footnote{For a system with two space dimensions, such as those for shallow water.}  preserving diffeomorphisms in the independent variables. 
The preservation of potential vorticity in a numerical scheme is desired but thought to be an impossible goal. 

In recent work with T.\ Pryer, it was shown that if the numerical variational method, specifically a Finite Element method, has the same set of Lie group symmetries, then there will be what 
we call here \textit{Noether identities}, which are the analogues of the identities yielding Noether's conservation laws in the regular case \cite{Mpryer}. 
Further, the Noether identities hold for any mesh, even when the Lie group action is on the base space.

We begin by discussing Noether's conservation laws and the conservation of potential vorticity in smooth (regular) variational systems. The aim is to provide essentially a tutorial in how Noether's identity is calculated, from the infinitesimal criterion of invariance,
in simple examples for readers unfamiliar with the, by now, extensive literature.

We then review the results that derive by applying this criterion to two families of variational numerical schemes.
 
Specifically, we  show first that the Finite Element Noether identity corresponding to invariance under the {particle relabelling} symmetry, does indeed converge to the regular law giving the conservation 
of potential vorticity, at the same time as having discrete analogues of energy and linear and angular momenta.  We then construct a finite difference type scheme which also demonstrates finite difference analogues
of these laws.

\section{Noether's laws for smooth variational systems}

A Lie group is a group which is also a manifold, a topological space with an atlas of smooth Euclidean coordinate charts. In practice, this means the group elements depend smoothly on real parameters. 
The two simplest examples occurring in physical and engineering problems are: \begin{enumerate}
\item  the real numbers under addition, $(\mathbb{R}, +)$, acting on
a variable $x$ by either translation, $x\mapsto x+\epsilon$,  or scaling, $x\mapsto \exp(\lambda ) x$
\item  the special orthogonal group, $SO(n)$, acting on a point $\mathbf{v}\in \mathbb{R}^n$ linearly, as $\mathbf{v} \mapsto g \mathbf{v}$, that is, as a rotation.
\end{enumerate} 

Any action on variables induces an action on functions of those variables, $g\cdot F(x)= F(g\cdot x)$. The induced action on derivatives is obtained via the chain rule, 
$$g\cdot \frac{\partial u}{\partial x}:=\frac{\partial \left(g\cdot u\right)}{\partial \left(g\cdot x\right)}=\frac{\ds\frac{\partial g\cdot u}{\partial  x}}{\ds\frac{\partial g\cdot x}{\partial x}},$$ 
where the second equality follows from implicit differentiation in the case here of a single independent variable.
Further, the induced action on integrands of line, surface and volume forms, is also induced by the chain rule, for example, $$g\cdot {\rm d}x = {\rm d} \left(g\cdot x\right)=\frac{\partial \left(g\cdot x\right)}{\partial  x}\,{\rm d}x.$$

\begin{example}\label{exone} Consider the action of $SE(2)=S)(2)\ltimes \mathbb{R}^2$, given by rotation and translation in the plane. If we take a curve in the plane to be $(x,u(x))$, then we have
$$g\cdot (x,u(x)) = \left(g\cdot x, g\cdot u\right)=\left(\cos(\theta)\, x - \sin(\theta)\, u(x) + a, \sin(\theta)\, x + \cos(\theta)\, u(x) +b\right)$$
The induced action on $u_x$ is
$$g\cdot u_x = \frac{{\rm d} (g\cdot u)}{{\rm d}( g\cdot x)} = \frac{ {\rm d}( g\cdot u)}{{\rm d} x}\Big/ \frac{{\rm d} (g\cdot x)}{{\rm d} x} = \frac{\sin(\theta) + \cos(\theta) u_x}{\cos(\theta) - \sin(\theta) u_x},$$
while that on $u_{xx}$ is 
$$g\cdot u_{xx} = \frac{ {\rm d}( g\cdot u_x)}{{\rm d} x}\Big/ \frac{{\rm d} (g\cdot x)}{{\rm d} x} = \frac{u_{xx}}{\left(\cos(\theta) - \sin(\theta)\, u_x\right)^3}.$$
The induced action on the integrand, $ f(x,u(x))\,{\rm d}x $ is given by
$$ \begin{array}{rcl} g\cdot \left( F(x,u(x))\,{\rm d}x\right)  &=& F(g\cdot x, g\cdot u)\, {\rm d}(g\cdot x)  \\
&=& F(g\cdot x, g\cdot u)\,\left( \cos(\theta) -\sin(\theta)\,  u_x\right){\rm d}x.\end{array}$$
\end{example}\smallskip

We are interested in the induced action on a Lagrangian, which for a one dimensional problem, takes the form,
\begin{equation}\label{actLag}g\cdot \int_{\Omega} L(x,u,u_x)\,{\rm d}x :=  \int_{g\cdot \Omega}  L\left(g\cdot x,g\cdot u,g\cdot u_x\right)\, {\rm d}\left(g\cdot x\right) = \int_{\Omega} L\left(g\cdot x,g\cdot u,g\cdot u_x\right)
\frac{\partial g\cdot x}{\partial  x}\,{\rm d}x\end{equation}
where the second equality follows by considering all quantities with respect to the original independent variable
This is the same as transforming back to the original variables on the base space (\cite{olver}, page 249) via the group action map on the base space. 
It is important to note that as $g\cdot x$ is a function of $(x, u(x),g)$, then the third integrand is indeed defined on the original domain $\Omega$.

If $g\cdot \int L\,{\rm d}x = \int L\,{\rm d}x$ for all $g\in G$, then we may take a path $\epsilon\mapsto \gamma(\epsilon)\subset G$ though the identity of the group, so that
$\gamma(0)=e$, the identity of $G$, and we will have
\begin{equation}\label{infcrit} \frac{{\rm d}}{{\rm d}\epsilon}\Big\vert_{\epsilon =0} \gamma(\epsilon) \cdot \int L\,{\rm d}x = 0.\end{equation}
The result is the so-called \textit{Infinitesimal Criterion of Invariance}. As we show below, Noether's Theorem is the straightforward result of  this criterion for an invariant Lagrangian.

\begin{example} Continuing with Example \ref{exone}, we have $g=g(\theta, a, b)$, and $g(0,0,0)=e$, the identity map.  By the formulae for the group action, we obtain the so-called infinitesimals
$$\frac{{\rm d}}{{\rm d} a }\Big\vert_{a =0} g\cdot x = 1,\qquad  \frac{{\rm d}}{{\rm d}a}\Big\vert_{a =0} g\cdot u =0,\qquad  \frac{{\rm d}}{{\rm d} a }\Big\vert_{a =0} g\cdot u_x = 0,\qquad \frac{{\rm d}}{{\rm d} a }\Big\vert_{a =0} g\cdot u_{xx} = 0 $$
and similarly for $b$ and $\theta$. We collect the infinitesimals into a table, where for variable $v$ and parameter $a_i$ the $(a_i, v)$ entry is $ {\rm d} v/{\rm d} a_i \vert_{e}$.
$$\begin{array}{|c|ccccc|}\hline
& x&u&u_x&u_{xx}& {\rm d}x\\ \hline 
a& 1&0&0&0&0\\
b& 0&1&0&0&0\\
\theta& -u & x & 1+u_x^2  &-3 u_x u_{xx}&  -u_x \,{\rm d}x\\ \hline
\end{array}$$
We then have,
$$ \frac{{\rm d}}{{\rm d}\theta}\Big\vert_{\theta =0} F\left(g\cdot x, g\cdot u \right) = - \frac{\partial F}{\partial x} u + \frac{\partial F}{\partial u} x$$
and 
$$ \frac{{\rm d}}{{\rm d}\theta}\Big\vert_{\theta =0} F\left(g\cdot x, g\cdot u \right)\,{\rm d}(g\cdot x) = \left(- \frac{\partial F}{\partial x} u + \frac{\partial F}{\partial u} x -F(x,u) u_x \right)\,{\rm d}x.$$
If we have a   second order Lagrangian which is invariant under this action, so that, $L\left(g\cdot x, g\cdot u, g\cdot u_x, g\cdot u_{xx} \right)\,{\rm d}(g\cdot x) = L(x,u,u_x,u_{xx})\,{\rm d}x$ for all $g$, then there are
three identities that $L$  satisfies,
$$\begin{array}{crcl}a: & \ds\frac{\partial L}{\partial x}  &=&0\\[7.5pt]
b: &\ds\frac{\partial L}{\partial u}  &=&0\\[7.5pt]
\theta: &\ds -u  \frac{\partial L}{\partial x} + x \frac{\partial L}{\partial u} + (1+u_x^2) \frac{\partial L}{\partial u_x} -3 u_x u_{xx} \frac{\partial L}{\partial u_{xx}} -u_x  L &=&0
\end{array}$$
\end{example}\bigskip

We now turn to the  statement of Noether's Theorem. \textbf{Noether's theorem is, in fact, ``nothing but" a rearrangement of the infinitesimal criterion of invariance.} 
We show only the simplest low order cases, and then those that pertain to our systems of interest. The text by Olver, gives a comprehensive description and proof of recurrence formulae for the infinitesimals, and the proof of the Noether's Theorem, in the general case. We record here the result for some
simple examples.

 Consider the simplest, lowest order, lowest dimension, Lagrangian, $L(x,u,u_x)\,{\rm d}x$,
and suppose this Lagrangian is invariant under a Lie group action with parameter $\epsilon$. If the infinitesimals are
\begin{equation}\label{InfDefsXiPhi}\frac{{\rm d}}{{\rm d} \epsilon }\Big\vert_{\epsilon =0} g(\epsilon) \cdot x = \xi(x,u),\qquad \frac{{\rm d}}{{\rm d} \epsilon }\Big\vert_{\epsilon =0} g(\epsilon) \cdot u = \phi(x,u)\end{equation}
then the infinitesimal criterion for invariance yields,
\begin{equation}\label{lowest}  \frac{\partial L}{\partial x}\xi(x,u)  + \frac{\partial L}{\partial u} \phi(x,u)+  \frac{\partial L}{\partial u_x}\left(\frac{{\rm d}\phi}{{\rm d}x} -u_x \frac{{\rm d}\xi}{{\rm d}x} \right) +  L \frac{{\rm d}\xi}{{\rm d}x} =0.\end{equation}

 Adding and subtracting the terms needed to convert $\partial L / \partial x  $ into ${\rm d}L/{\rm d}x$ in the Infinitesimal Criterion of Invariance, Equation (\ref{lowest}) and  rearranging, we obtain
\begin{equation}\label{NoetherSmooth} \left(\phi -u_x \xi\right)E(L) + \frac{{\rm d}}{{\rm d}x}\left(L\xi + \frac{\partial L}{\partial u_x} \left(\phi -u_x \xi\right)\right) =0\end{equation}
where $$E(L)=\frac{\partial L}{\partial u}-
\frac{{\rm d}}{{\rm d}x} \frac{\partial L}{\partial u_x}$$ is the Euler Lagrange equation. Therefore, on solutions of $E(L)=0$, we will have
$$\frac{{\rm d}}{{\rm d}x}\left(L\xi + \frac{\partial L}{\partial u_x} \left(\phi -u_x \xi\right)\right) =0$$
which is Noether's conservation law for this case.

\begin{example} For a second order, one dimensional Lagrangian, $L(x,u,u_x,u_{xx})\,{\rm d}x$, which is invariant under a Lie group action with infinitesimals $\xi$, $\phi$, Noether's Theorem is, setting $Q=\phi-u_x\xi$,
\begin{equation}\label{casexyuxx}E(L)Q +\frac{{\rm d}}{{\rm d}x}\left( L\xi + \frac{\partial L}{\partial u_x} Q + \frac{\partial L}{\partial u_{xx}} \frac{{\rm d}}{{\rm d}x}Q-Q \frac{{\rm d}}{{\rm d}x}\frac{\partial L}{\partial u_{xx}}    \right)=0,\end{equation}
where $$E(L)=\frac{\partial L}{\partial u}-
\frac{{\rm d}}{{\rm d}x} \frac{\partial L}{\partial u_x}+\frac{{\rm d}^2}{{\rm d}x^2} \frac{\partial L}{\partial u_{xx}}$$ is the Euler Lagrange equation for this case.
Equation (\ref{casexyuxx}) is easily seen to be a rearrangement of the infinitesimal criterion of invariance, which is
$$\frac{\partial L}{\partial x}\xi +\frac{\partial L}{\partial u}\phi +\frac{\partial L}{\partial u_x}\phi_{[x]} +\frac{\partial L}{\partial u_{xx}}\phi_{[xx]} +L\frac{{\rm d}}{{\rm d} x}\xi=0$$
where 
$$\begin{array}{rclcl}
\phi_{[x]}&=&\ds\frac{{\rm d}}{{\rm d}\epsilon}\big\vert_{\epsilon=0} \widetilde{u_x}  &=& \ds\frac{{\rm d}\phi}{{\rm d} x}-u_x \frac{{\rm d}\xi}{{\rm d} x}\\[7.5pt]
\phi_{[xx]}&=&\ds\frac{{\rm d}}{{\rm d}\epsilon}\big\vert_{\epsilon=0} \widetilde{u_{xx}}  &=&\ds\frac{{\rm d}}{{\rm d}x}\left( \frac{{\rm d}\phi}{{\rm d} x}-u_x \frac{{\rm d}\xi}{{\rm d} x}\right) -u_{xx}\frac{{\rm d}\xi}{{\rm d} x}
\end{array}$$
\end{example}

\begin{remark}[Names of the conserved quantities]
If the independent variables are $(t,x_1,x_2,\dots, x_p)$ where $t$ is time, and $u^{\alpha}$, $\alpha=1,\dots, q$ are the dependent variables, then
Noether's law is of the form $$\sum_{\alpha}Q^{\alpha}E^{\alpha}(L)+ \frac{{\rm d}}{{\rm d}t} A_0+\sum_i \frac{{\rm d}}{{\rm d}x_i}  A_i=0,$$
where $E^{\alpha}(L)$ is the Euler Lagrange equation with respect to the dependent variable $u^{\alpha}$ and $Q^{\alpha}$ is the characteristic of the group action with respect to $u^{\alpha}$.
If  $\Omega\subset \mathbb{R}^p$ is the spatial domain for the Lagrangian, we have on solutions of the Euler Lagrnage equations,
$$0= \int_{\Omega}\left( \frac{{\rm d}}{{\rm d}t} A_0+\sum_i \frac{{\rm d}}{{\rm d}x_i}  A_i \right)  
= \frac{{\rm d}}{{\rm d}t} \int_{\Omega}A_0 + \int_{\partial \Omega} \mathbf{A}\cdot\mathbf{n}$$
by Stokes' Theorem. The scalar $A_0$ is said to be \textit{conserved density} while the vector $(A_1, \dots, A_p)$ is called the \textit{flux}. 
If the Lie group action is translation in time, the scalar $A_0$ is said to be the \textit{energy} of the system described by the Euler Lagrange equations.
The table following gives the common names for the conservation laws for the given symmetries.
\newline
 \begin{center}
\smaller{
\begin{tabular}{|ll|}
\hline
{\bf Symmetry}& {\bf $A_0$}\\
Translation in time & Energy\\
Translation in space & Linear momentum\\
Rotation in space & Angular momentum\\
Particle relabelling & Potential vorticity\\ \hline
\end{tabular}}\end{center}
\end{remark}

\subsection{Shallow water Lagrangians}
We are concerned with Lagrangians for shallow water systems. The general form of the infinitesimal criterion of invariance and Noether's laws for such systems  are detailed in the next example.
\begin{example}\label{SWWone}[Lagrangians of shallow water type - general form]  We consider the case in which there are two spatial and one time variables, which we will denote as $(a,b,t)$ and which will be the particle labels and time. We assume two dependent variables, which we will denote $(x,y)$ and which will be the position of the fluid particle at time $t$, so that $(x(a,b,0),y(a,b,0))=(a,b)$ and we will assume a first order Lagrangian, that is, 
\begin{equation}\label{LagOrder1simple}L(a,b,t,x_a,x_b,x_t,y_a,y_b,y_t)\, {\rm d}a\,{\rm d}b\,{\rm d}t.\end{equation} We take the infinitesimals to be
\begin{equation}\label{InfXiPhi}\frac{{\rm d}}{{\rm d} \epsilon}\big\vert_{\epsilon=0} \widetilde{a}= \xi, \quad \frac{{\rm d}}{{\rm d} \epsilon}\big\vert_{\epsilon=0} \widetilde{b}= \tau,\quad \frac{{\rm d}}{{\rm d} \epsilon}\big\vert_{\epsilon=0} \widetilde{t}= \chi,\quad  \frac{{\rm d}}{{\rm d} \epsilon}\big\vert_{\epsilon=0} \widetilde{x}= \phi^x,\quad  \frac{{\rm d}}{{\rm d} \epsilon}\big\vert_{\epsilon=0} \widetilde{y}= \phi^y.\end{equation}
Then the infinitesimal criterion of invariance is
$$0=\frac{\partial L}{\partial a}\xi +\frac{\partial L}{\partial b}\tau+\frac{\partial L}{\partial t}\chi+\frac{\partial L}{\partial x}\phi^x + \frac{\partial L}{\partial y}\phi^y+\frac{\partial L}{\partial x_a}\phi^x_{[a]}+
\frac{\partial L}{\partial x_b}\phi^x_{[b]}+\frac{\partial L}{\partial x_t}\phi^x_{[t]}+\frac{\partial L}{\partial y_a}\phi^y_{[a]}+\frac{\partial L}{\partial y_b}\phi^y_{[b]}+\frac{\partial L}{\partial y_t}\phi^y_{[t]}$$
where
$$\begin{array}{rclcl}\phi^x_{[a]}&=&\ds \frac{{\rm d}}{{\rm d}\epsilon}\big\vert_{\epsilon=0} \widetilde{x_a}&=&\ds \frac{{\rm d}}{{\rm d} a}\phi^x - x_a\frac{{\rm d}}{{\rm d} a}\xi -x_b \frac{{\rm d}}{{\rm d} a}\tau-x_t\frac{{\rm d}}{{\rm d} a}\chi\\[7.5pt]
\phi^x_{[b]}&=&\ds \frac{{\rm d}}{{\rm d}\epsilon}\big\vert_{\epsilon=0} \widetilde{x_b}&=&\ds \frac{{\rm d}}{{\rm d} b}\phi^x - x_a\frac{{\rm d}}{{\rm d} b}\xi -x_b \frac{{\rm d}}{{\rm d} b}\tau-x_t\frac{{\rm d}}{{\rm d} b}\chi\\[7.5pt]
\phi^x_{[t]}&=&\ds \frac{{\rm d}}{{\rm d}\epsilon}\big\vert_{\epsilon=0} \widetilde{x_t}&=&\ds \frac{{\rm d}}{{\rm d} t}\phi^x - x_a\frac{{\rm d}}{{\rm d} t}\xi -x_b \frac{{\rm d}}{{\rm d} t}\tau-x_t\frac{{\rm d}}{{\rm d} t}\chi\\[7.5pt]
\phi^y_{[a]}&=&\ds \frac{{\rm d}}{{\rm d}\epsilon}\big\vert_{\epsilon=0} \widetilde{y_a}&=&\ds \frac{{\rm d}}{{\rm d} a}\phi^y - y_a\frac{{\rm d}}{{\rm d} a}\xi -y_b \frac{{\rm d}}{{\rm d} a}\tau-y_t\frac{{\rm d}}{{\rm d} a}\chi\\[7.5pt]
\phi^y_{[b]}&=&\ds \frac{{\rm d}}{{\rm d}\epsilon}\big\vert_{\epsilon=0} \widetilde{y_b}&=&\ds \frac{{\rm d}}{{\rm d} b}\phi^y - y_a\frac{{\rm d}}{{\rm d} b}\xi -y_b \frac{{\rm d}}{{\rm d} b}\tau-y_t\frac{{\rm d}}{{\rm d} b}\chi\\[7.5pt]
\phi^y_{[t]}&=&\ds \frac{{\rm d}}{{\rm d}\epsilon}\big\vert_{\epsilon=0} \widetilde{y_t}&=&\ds \frac{{\rm d}}{{\rm d} t}\phi^y - y_a\frac{{\rm d}}{{\rm d} t}\xi -y_b \frac{{\rm d}}{{\rm d} t}\tau-y_t\frac{{\rm d}}{{\rm d} t}\chi\end{array}
$$
The infinitesimal criterion is then rearranged to yield Noether's identity for this case, which is
\begin{align}\label{caseabtxy}
E^x(L)Q^x &+ E^y(L)Q^y +\ds\frac{{\rm d}}{{\rm d}a}\left( L\xi +\frac{\partial L}{\partial x_a}Q^x 
+\frac{\partial L}{\partial y_a}Q^y \right)\nonumber\\ 
&+\frac{{\rm d}}{{\rm d}b}\left(L\tau +\frac{\partial L}{\partial x_b}Q^x+\frac{\partial L}{\partial y_b}Q^y  \right) +\frac{{\rm d}}{{\rm d}t}\left(L\chi +\frac{\partial L}{\partial x_t}Q^x+\frac{\partial L}{\partial y_t}Q^y  \right) =0
\end{align}
where   $Q^x=\phi^x - x_a\xi-x_b\tau-x_t\chi$,  $Q^y=\phi^y - y_a\xi-y_b\tau-y_t\chi$ and where the Euler Lagrange equations are
$$\begin{array}{rcl} E^x(L)&=&\ds  \ds\frac{\partial L}{\partial x} - \frac{{\rm d}}{{\rm d}a}\frac{\partial L}{\partial x_a} - \frac{{\rm d}}{{\rm d}b}\frac{\partial L}{\partial x_b} - \frac{{\rm d}}{{\rm d}t}\frac{\partial L}{\partial x_t} \\[12pt]
E^y(L)&=&\ds \ds \frac{\partial L}{\partial y} - \frac{{\rm d}}{{\rm d}a}\frac{\partial L}{\partial y_a} - \frac{{\rm d}}{{\rm d}b}\frac{\partial L}{\partial y_b} - \frac{{\rm d}}{{\rm d}t}\frac{\partial L}{\partial y_t}.
\end{array}$$
\end{example}




We next  look
at the particle relabelling symmetry and the conservation of potential vorticity for shallow water Lagrangians.

\subsubsection{The particle relabelling symmetry and the conservation of potential vorticity}
We consider again the problem of Example \ref{SWWone}, with independent variables $(a,b,t)$ where $t$ is time and $(a,b)$ are the particle labels at time $t=0$. 
The dependent variables $(x,y)$ are the position of the particle at time $t$. The particle relabelling symmetry group is given by the set of area preserving diffeomorphisms,
$$\widetilde{a} = A(a,b),\qquad \widetilde{b}=B(a,b), \qquad A_aB_b-A_bB_a=1.$$
Further, $\widetilde{t}=t$, $\widetilde{x}=x$ and $\widetilde{y}=y$, that is, $x$, $y$ and $t$ are invariant under the action. 
A first order Lagrangian which is invariant under this
symmetry action takes the form
\begin{equation}\label{SWLagwithPotV} L(x,y,\Delta, x_t,y_t)\,{\rm d}a\,{\rm d}b\,{\rm d}t,\qquad \Delta = x_ay_b-x_by_a.\end{equation}
A higher order Lagrangian with this symmetry will have additional arguments of the form of invariant derivatives of the generating invariant $\Delta$,
\begin{equation}\label{invDiffOps}
\frac{\partial^{m+n}}{\partial x^m\partial y^n} \Delta,\qquad \frac{\partial}{\partial x} =\frac1{\Delta}\left( y_b\frac{\partial}{\partial a}-y_a\frac{\partial}{\partial b}\right),\quad \frac{\partial}{\partial y} =\frac1{\Delta}\left( -x_b\frac{\partial}{\partial a}+x_a\frac{\partial}{\partial b}\right)\end{equation}
as well as derivatives with respect to $t$. It should noted that $\ds\left[\frac{\partial}{\partial t},\frac{\partial}{\partial x}\right]$ and $\ds\left[\frac{\partial}{\partial t},\frac{\partial}{\partial y}\right]$ are not zero, as $x$ and $y$ depend on time $t$.

The formula for the  potential vorticity of the smooth problem is obtained by cross-differentiation of the formulae for the linear momenta, an observation by Peter Hydon, \cite{hydon}. 
The Noether identities for the smooth linear momenta of a first order Lagrangian, are given by
$$\begin{array}{rcl}
x_a E^x + y_aE^y &=&\ds  \frac{{\rm d}}{{\rm d}a} \left(L+\Delta \frac{\partial L}{\partial \Delta}\right) + 
\frac{{\rm d}}{{\rm d}t}\left( x_a\frac{\partial L}{\partial x_t} + y_a\frac{\partial L}{\partial y_t}\right)\\[7.5pt]
x_b E^x + y_bE^y &=&\ds  \frac{{\rm d}}{{\rm d}b} \left(L+\Delta \frac{\partial L}{\partial \Delta}\right) + 
\frac{{\rm d}}{{\rm d}t}\left( x_b\frac{\partial L}{\partial x_t} + y_b\frac{\partial L}{\partial y_t}\right).\end{array}
$$ Cross differentiation of these, on solution of the Euler Lagrange equations, yields the conservation of potential vorticity, given explicitly as
\begin{equation}\label{potvortsmoothone}\frac{{\rm d}}{{\rm d}t}\left( \frac{{\rm d}}{{\rm d}b} \left( x_a\frac{\partial L}{\partial x_t} + y_a\frac{\partial L}{\partial y_t}\right) -\frac{{\rm d}}{{\rm d}a}\left( x_b\frac{\partial L}{\partial x_t} + y_b\frac{\partial L}{\partial y_t}\right) \right) =0,\end{equation}
which can be simplified to 
\begin{equation}\label{potvortsmoothoneTWO}\frac{{\rm d}}{{\rm d}t}\left( \frac{\partial}{\partial y} \frac{\partial L}{\partial x_t} -\frac{\partial }{\partial x}\frac{\partial L}{\partial y_t} \right)=0\end{equation}
using the definitions of the operators in Equation (\ref{invDiffOps}).
The exactly conserved density in Equation (\ref{potvortsmoothoneTWO}) is known as potential vorticity $\Psi$,
$$\Psi= \frac{\partial}{\partial y} \frac{\partial L}{\partial x_t} -\frac{\partial }{\partial x}\frac{\partial L}{\partial y_t} $$
We next write down two well known shallow water Lagrangians and calculate their potential vorticities.

\subsection{Examples}\label{mainExssubsec}
We give here two examples, quoted from \cite{BilaMC}.

 In the
two-dimensional shallow water theory \cite{refRubRoul} a typical particle
(more precisely, fluid column) has the Cartesian horizontal coordinates
\begin{equation} x=x(a,b,t),\qquad y=y(a,b,t),
\label{eq:1.1}\end{equation} expressed as functions of the
particle labels $(a,b)\in \R^2 $ and time $t\in \R^{+}$. For
convenience, each particle is labelled by its position at a reference
time $t=0$, this means the functions in (\ref{eq:1.1}) are defined such
that $x(a,b,0)=a$ and $y(a,b,0)=b$. The incompressibility hypothesis
requires
\begin{equation}{\frac{h(a,b,0)}{h(a,b,t)}}={\frac{\partial (x,y)}{\partial
(a,b)}},\label{eq:1.2}\end{equation}
where the Jacobian on the right is that of the mapping (\ref{eq:1.1}).
The time derivative of (\ref{eq:1.2}) following the particle gives the
continuity equation. In this paper we assume $h(a,b,0)=1$, so that the
incompressibility hypothesis becomes
\begin{equation}
h(a,b,t)=\frac{1}{x_ay_b-x_by_a},
\label{eq:1.3}\end{equation}
where subscripts denote partial derivatives. The mapping (\ref{eq:1.1}) is
assumed to be invertible, such that when $a=a(x,y,t)$ and $b=b(x,y,t)$ are
inserted into (\ref{eq:1.2}), then the current depth $h$ is expressed as
function of $x$, $y$ and $t$, which represents the Eulerian description.

\noindent\textbf{Example 1: the potential form of the \SW}

The equations of the horizontal momentum balance for the flows
over a bed which is rotating with position dependent Coriolis
parameter $f=f(y)$ are
\begin{subequations}\begin{align}
\ddot x +gh_x-f\dot{y}&=0,\\
\ddot y +gh_y+f\dot{x}&=0,
\end{align}
\label{eq:1.4}\end{subequations}
where $g$ is a nonzero constant (representing the combined effect of the
acceleration due to gravity and a centrifugal component due to the Earth's
rotation), dot denotes the time derivative following a particle, and $h_x$
and $h_y$ are given by
\begin{subequations}\begin{align}
h_x&=h(y_bh_a-y_ah_b),\\
h_y&=h(x_ah_b-x_bh_a).
\end{align}
\label{eq:1.5}\end{subequations}
Henceforth we shall assume that $f$ is a constant. It is known that
the {\it shallow water potential vorticity} defined by
\begin{equation}
\Omega ={\frac1h}\left (\frac{\partial \dot {y}}{\partial
x}-\frac{\partial \dot {x}}{\partial y }+f\right),
\label{eq:1.6}\end{equation}
is conserved on particles (see \cite{refBHR,refPM} and the references
therein).

If we let $u=\dot{x}$ and $v=\dot{y}$, then the \SW\ (\ref{eq:1.4}) can be
written as
\begin{subequations}\begin{align}
\dot {x} &= u,\\
\dot {y}&= v,\\
\dot {u} &+gh(y_bh_a-y_ah_b)-fv=0,\\
\dot {v} &+gh(x_ah_b-x_bh_a)+fu=0,
\end{align}\label{eq:3.11}\end{subequations}
where $h$ is given by (\ref{eq:1.3}), which we shall refer to as the {\it
potential form} of the \SW\ (\ref{eq:1.4}).

Salmon \cite{refSalmon} showed that there is a first-order Lagrangian
\begin{equation}
\L=(u-R)\dot{x}+(v+P)\dot{y}-\tfrac12 \left(u^2+v^2+gh\right),
\label{eq:6.1}\end{equation}
where the functions $P=P(x,y)$ and $R=R(x,y)$ satisfy
$$P_x+R_y=f, $$
with $h$ given by (\ref{eq:1.3}) and $f$ constant, for which the associated
Euler-Lagrange equations represent (up to a sign) the \SW\
(\ref{eq:3.11}), that is, the potential form of the \SW.

It can be checked that the formula for the conserved quantity is indeed that given in Equation (\ref{eq:1.6}), as the new invariants $u$ and $v$ make no essential difference to the computation.

\noindent\textbf{Example 2: the semi-geostrophic \SW}
In the past thirty years, the semi-geostrophic equations have become a
model for describing atmospheric motions on a synoptic scale, including
the presence of fronts \cite{refRoulNorb,refRoulSew97}.

The semi-geostrophic approximation to equations (\ref{eq:1.4}) is the
replacement of the true acceleration by the time derivative of the
vector
\begin{equation}
u_g=-\frac{g}{f}h_y,\qquad
v_g=\frac{g}{f}h_x, 
\label{eq:1.7}\end{equation}
following the particle. This vector field 
is called the {\it geostrophic velocity}. Thus, the
semi-geostrophic approximation seeks to find motions satisfying
\begin{subequations}\begin{align}
\dot{u}_g+gh_x-f\dot{y}&=0,\\
\dot{v}_g+gh_y+f\dot{x}&=0.
\end{align}\label{eq:1.8}\end{subequations}
It is also known that the {\it semi-geostrophic potential
vorticity}, given by
\begin{equation}
\Omega^{*}=\frac{1}{h}\left (f+\frac{\partial v_g}{\partial
x}-\frac{\partial u_g}{\partial y }+\frac1{f} \frac{\partial
(u_g,v_g)}{\partial (x,y)}\right),
\label{eq:1.9}\end{equation}
is conserved on particles \cite{refRoulSew96}.

The {\it classical form} of the \SG\ (\ref{eq:1.8}) are
\begin{subequations}\begin{align} hh_b\dot {x_a}-hh_a\dot {x_b}+(\dot
hh_b+h\dot {h_b})x_a-(\dot hh_a+h\dot {h_a})x_b+fhh_by_a-fhh_ay_b&=-
\frac{f^2}{g}\dot{y},\\ hh_b\dot {y_a}-hh_a\dot {y_b}+(\dot hh_b+h\dot
{h_b})y_a-(\dot hh_a+h\dot {h_a})y_b-fhh_bx_a +fhh_ax_b&=\frac{f^2}
{g}\dot{x}, \end{align}\label{eq:4.1}\end{subequations} with
$h$ given by (\ref{eq:1.3}). This system is obtained after substituting
the functions $u_g$ and $v_g$ given by (\ref{eq:1.7}) into (\ref{eq:1.8}).

The potential form of the \SG\ is given by 
(\ref{eq:1.7}) and (\ref{eq:1.8}), i.e.
\begin{subequations}\begin{align}
u_g&=-\frac{g}{f}h(x_ah_b-x_bh_a),\\
v_g&=\frac{g}{f}h(y_bh_a-y_ah_b),\\
\dot{u}_g&+gh(y_bh_a-y_ah_b)-f\dot{y}=0,\\
\dot{v}_g&+gh(x_ah_b-x_bh_a)+f\dot{x}=0,
\end{align}\label{eq:4.6}\end{subequations}
where $f$ is a nonzero constant and $h$ is given by (\ref{eq:1.3}).

The
\SW\ (\ref{eq:4.6}) represent the Euler-Lagrange equations
associated with the first-order Lagrangian
\cite{refHydon1}\begin{equation}
\L=(u_g-R)\dot {x}+(v_g+P)\dot {y}-\tfrac12(u_g^2+v_g^2+gh)-r\dot
{u}_g+p\dot {v}_g, \label{eq:7.1}\end{equation}
where $P(x,y)$, $R(x,y)$, $p(u_g,v_g)$ and $r(u_g,v_g)$ are arbitrary
functions satisfying
\begin{equation}P_x+R_y=f,\qquad
p_{u_g}+r_{v_g}={1/f},\label{eq:7.2}\end{equation}
and $h$ is given by (\ref{eq:1.3}).

This case requires an extension of the theory elaborated above as we have additional invariants with time derivatives. Working through the calculations of the conserved quantity, we obtain
\begin{equation}\label{potvortsmoothoneSGPOT}\frac{{\rm d}}{{\rm d}t}\left( \frac{{\rm d}}{{\rm d}b} \left( x_a\frac{\partial L}{\partial x_t} + y_a\frac{\partial L}{\partial y_t}  + u_{g,a}\frac{\partial L}{\partial \dot{u_g}}+
v_{g,a}\frac{\partial L}{\partial \dot{v_g}}\right) 
-\frac{{\rm d}}{{\rm d}a}\left( x_b\frac{\partial L}{\partial x_t} + y_b\frac{\partial L}{\partial y_t}  + u_{g,b}\frac{\partial L}{\partial \dot{u_g}}+
v_{g,b}\frac{\partial L}{\partial \dot{v_g}}\right) \right) =0,\end{equation}
The conserved quantity is straightforwardly shown to be that given in Equation \ref{eq:1.9}, using the definitions of the operators in Equation (\ref{invDiffOps}) and noting that, for example,
$$ \frac{\partial}{\partial x}\frac{\partial L}{\partial \dot{u_g}} =  \frac{\partial r}{\partial x} = \frac{\partial r}{\partial u_g}\frac{\partial u_g}{\partial x}+\frac{\partial r}{\partial v_g}\frac{\partial v_g}{\partial x}$$
and similarly for 
$$ \frac{\partial}{\partial y}\frac{\partial L}{\partial \dot{u_g}}, \qquad \frac{\partial}{\partial x}\frac{\partial L}{\partial \dot{v_g}},\qquad \frac{\partial}{\partial y}\frac{\partial L}{\partial \dot{v_g}}$$
appearing in Equation (\ref{potvortsmoothoneSGPOT}), and finally using Equation (\ref{eq:7.2}).

\section{Noether identities for numerical variational problems -- Finite Element}

We now look at the Finite Element Noether identity, developed by the author with T.\ Pryer, \cite{Mpryer}. We treat the approximate Lagrangian as an exact problem, and require invariance of the Finite Element Lagrangian. 
We no longer have point wise invariance, and so the Finite Element Noether identity will be an integral expression.

We consider the Lagrangian, $\int L(x,u,u_x)\, {\rm d}x$. The Finite Element Euler Lagrange equations for the Finite Element solution $U$ are
\begin{equation}\label{FEMEL}
0=\int_{\Omega} E^u\big\vert_{u=U} V +\int_{\mathcal{F}} \jump{
\frac{\partial L}{\partial U_x}} 
V, 
\qquad \forall V\in \mathbb{V}.
\end{equation} where $\mathbb{V}$ is an appropriate  space of piecewise continuous polynomial functions over a
triangulation, and where $E^u\vert_{u=U} = -\left({{\rm d}}/{{\rm d}x}\right){\partial L}/{\partial U_x}+{\partial L}/{\partial U}.$ 
The Finite Element Noether Identity is
\begin{equation}\label{FEMconlaw}
  \begin{split}
    0 
    &
= 
    \int_{\Omega} 
     E^u\big\vert_{u=U}
\left( \phi - U_x \xi\right)
    + \int_{\mathcal F}
    \jump{
      \frac{\partial L}{\partial U_x} 
      \left( \phi - U_x \xi \right) + L\xi
    }
  \end{split}
\end{equation}
where $\xi$ and $\phi$ are the infinitesimals defined in Equation (\ref{InfDefsXiPhi}).
 Thus, for example, for the conservation of linear momentum, where the group action is translation in $x$, so that $\widetilde{x}=x+\epsilon$, we must have that $L$ does not depend explicitly on $x$. Then the infinitesimals of Equation (\ref{InfXiPhi}) are given by
$(\xi,\phi)=(1,0)$  and the Noether identity for this symmetry is given by 
\begin{equation}\label{FEMconlawLM}
  \begin{split}
    0 
    &
= 
    \int_{\Omega} 
    E^u\big\vert_{u=U}
    \left( - U_x \right)
    + \int_{\mathcal F}
    \jump{
      \frac{\partial L}{\partial U_x} 
      \left(  - U_x  \right) + L
    }.
  \end{split}
\end{equation} This compares to the smooth case where we would have
${\rm d}/{\rm d}x \left(L-u_x \partial L/\partial u_x\right) =0$. We note that the term incorporating the Euler Lagrange equations in the FE  identity, cannot be set to zero, as it is not the case that $\phi-U_X\xi$ is in the test space in general.

We now turn to the FE case and consider our FE Noether identities given by the particle relabelling symmetry, for the Lagrangian in (\ref{SWLagwithPotV}). 
The infinitesimals needed to express the identity, analogous to those in Equation (\ref{InfXiPhi}), for $(a,b,t,x,y)$ are $(\xi^a, \xi^b,\xi^t,\phi^x,\phi^y)$,  and satisfy
$$\frac{\partial \xi^a}{\partial a} + \frac{\partial\xi^b}{\partial b}=0,\quad \xi^t=0,\quad \phi^x=0,\quad \phi^y=0.$$
We suppose  that the domain 
$\overline{\Omega}$ is made up of domains of the form $\overline{\Omega}_{i,n}=\Omega_i\times [t_n,t_{n+1}]$ where $\Omega_i\subset (a,b)$-space. 

\includegraphics[trim=100 100 200 0 0,clip=true,angle=-90,scale=0.5]{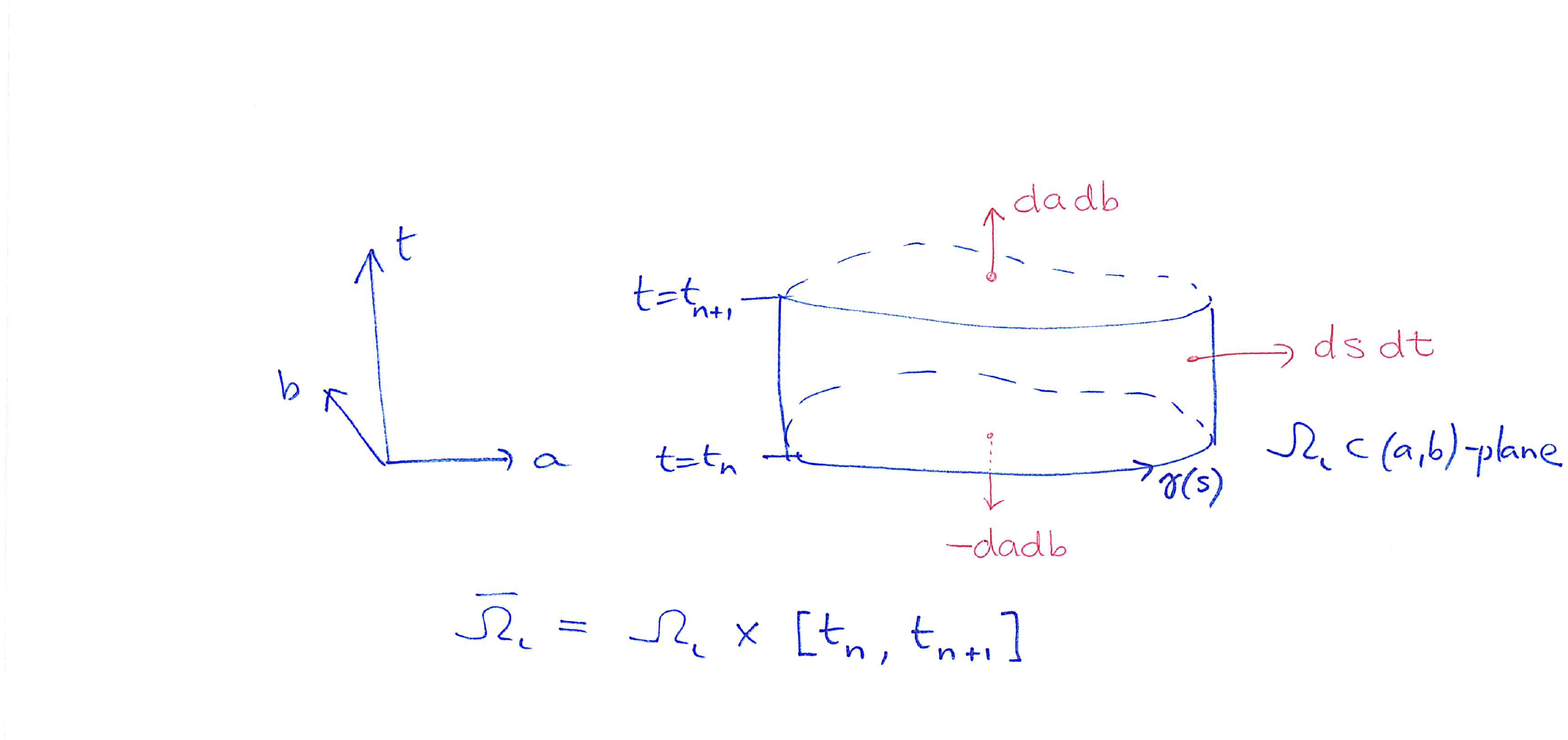}

Locally we may set 
$\xi^a=\partial_b \varphi=\varphi_b$ and $\xi^b=-\partial_a\varphi=\varphi_a$ for some smooth function $\varphi$, and then the FE Noether identity in this case, takes the form
\begin{equation}\label{mainFEPV}\begin{array}{rcl}
0&=&\ds \sum_i \int_{\overline{\Omega}_{i,n}} \left[\left(\varphi_aX_b-\varphi_bX_a\right)E^x\vert_{x=X} + \left(\varphi_aY_b-\varphi_bY_a\right)E^y\vert_{y=Y} \right]\, {\rm d}a\,{\rm d}b\,{\rm d}t\\[12pt]
&&\ds\qquad + \left(\int_{\Omega_i\times \{t_{n+1}\}} - \int_{\Omega_i\times \{t_{n}\}}\right)
\left(\left(\varphi_aX_b-\varphi_bX_a\right)\frac{\partial L}{\partial X_t} + 
\left(\varphi_aY_b-\varphi_bY_a\right) \frac{\partial L}{\partial Y_t}\right) \, {\rm d}a\,{\rm d}b\\[12pt]
&&\ds\qquad + \int_{t_n}^{t_{n+1}} \int_{\gamma_i(s)=\partial \Omega_i}
\left[\kern-2pt\left[ L+\Delta \frac{\partial L}{\partial \Delta}\right]\kern-2pt\right]\Big\vert_{\gamma_i} \frac{{\rm d}}{{\rm d}s} \varphi(\gamma_i(s))\, {\rm d}s\,{\rm d}t.
\end{array}\end{equation}
This is simply the three dimensional analogue of the above equation (\ref{FEMconlaw}), adapted to two dependent variables and the given infinitesimals for this group action.

We now show that the continuum limit of the identity (\ref{mainFEPV}) is indeed the conservation law (\ref{potvortsmoothoneTWO}).

If $|t_{n+1}-t_n|$ is small, we have for any continuous integrand $F$, that there exists $\tau$, $t_n<\tau<t_{n+1}$ such that
$$\int_{t_n}^{t_{n+1}} F(t)\,{\rm d}t = (t_{n+1}-t_n)F(\tau),$$
by the mean value theorem. We apply this result to the first and third integrals in Equation (\ref{mainFEPV}).
Doing this, then dividing through by $(t_{n+1}-t_n)$ and sending $t_{n+1}\rightarrow t_n$, the first and third intervals limit to the integral over $\Omega_i$ at $t=t_{n}$. The second integral in Equation (\ref{mainFEPV}) 
limits to a derivative with respect to $t$. We thus obtain, for any smooth function $\varphi$, that
\begin{equation}\label{mainFEPV2}\begin{split}
0&=\ds \sum_i \int_{\bar{\Omega}_{i}} \left[\left(\varphi_aX_b-\varphi_bX_a\right)E^x\vert_{x=X} + \left(\varphi_aY_b-\varphi_bY_a\right)E^y\vert_{y=Y} \right]\, {\rm d}a\,{\rm d}b\\
&\ds\qquad + \frac{{\rm d}}{{\rm d}t}\int_{\Omega_i} 
\left\{ \varphi_a\left(X_b \frac{\partial L}{\partial X_t}+Y_b\frac{\partial L}{\partial Y_t}\right)
-\varphi_b\left(X_a \frac{\partial L}{\partial X_t}+Y_a\frac{\partial L}{\partial Y_t}  \right)\right\} {\rm d}a\,{\rm d}b\\
&\ds\qquad +  \int_{\gamma_i(s)=\partial \Omega_i}
\left.\left[\kern-2pt\left[ L+\Delta \frac{\partial L}{\partial \Delta}\right]\kern-2pt\right]\right\vert_{\gamma_i} \frac{{\rm d}}{{\rm d}s} \varphi(\gamma_i(s))\, {\rm d}s.
\end{split}\end{equation}
Finally, we observe that Equation (\ref{mainFEPV2}) is the weak form of the smooth law. Indeed, for a strong solution of the Euler Lagrange equations, the $E^x$, $E^y$ terms are zero, as is the jump term in the third integrand. Taking $\varphi$ to be a bump function with support anywhere inside a mesh simplex, then 
integration by parts yields the same cross-derivative as that yielding the smooth potential vorticity, to be zero. 

Equation (\ref{mainFEPV2}) is straightforwardly extended to the potential forms of the equations and their Lagrangians in the examples in Section \ref{mainExssubsec}.

\def\refjl#1#2#3#4#5#6#7{\vspace{-0.25cm}
\bibitem{#1}{\frenchspacing#2}, {#3},
{\frenchspacing\it#4}, {\bf#5}\ (#7) #6.}

\def\refbk#1#2#3#4#5{\vspace{-0.25cm}
\bibitem{#1}{\frenchspacing\rm#2}, {\it #3\/}, #4 (#5).}

\def\refpp#1#2#3#4#5{\vspace{-0.25cm}
\bibitem{#1}{\frenchspacing\rm#2}, {\it #3\/}, #4 (#5).}

\def\refcf#1#2#3#4#5#6#7{\vspace{-0.25cm}
\bibitem{#1}{\frenchspacing\rm#2}, {\it #3\/},
in: ``{#4\/}"\ {\frenchspacing\rm#5}, #6 (#7).}

\end{document}